\documentclass[preprint,11pt]{elsarticle}
\usepackage{amssymb,amsmath,amsthm,amsfonts,amscd}
\usepackage{graphicx,color}
\usepackage{multicol}
\usepackage{float}
\usepackage{enumitem}
\usepackage{cancel} 
\journal{xxxxxxxxxxxxxxxxxx}
\newtheorem{theorem}{Theorem}[section]

\newtheorem{remark}{Remark}[section]

\begin{document}
	
	\begin{frontmatter}
		
		\title{Theoretical analysis of a multi-objective controllability problem for the linear wave equation on a non-cylindrical domain}
		\author[Pedro]{Pedro Paulo  A. Oliveira\corref{cor1}}
		\ead{pedropaulo@ufpi.edu.br}
		\author[Jesus]{Isa\'{i}as P.  de Jesus}
		\ead{isaias@ufpi.edu.br}
		\author[Neto]{Gilcenio R. de Sousa-Neto}
		\ead{gilceniorodrigues@gmail.com}
		\address[Pedro]{Federal University of Piau\'{\i}, DM, PI, Brazil}
		\address[Jesus]{Federal University of  Piau\'{\i}, DM, PI, Brazil}
		\address[Neto]{Federal University of Piau\'{\i}, DM, PI, Brazil}
		\cortext[cor1]{Corresponding author}
		
		\begin{abstract}
			In this article, we investigate certain theoretical aspects of the hierarchical controllability problem in one-dimensional wave equations within a moving domain using Stackelberg strategy. The controls are applied along a portion of the boundary and establish an equilibrium strategy among them, considering a leader control and a follower. We consider a linear wave equation.
		\end{abstract}
		
		
	\end{frontmatter}
	
	\section{Introduction}
	
	Let $Q_{T,\alpha}$ be a non-cylindrical domain associated with a function $\alpha>0$ and a constant $T>0$, defined as follows:
	$$
	Q_{T,\alpha} =
	\{
	(x,t)\in\mathbb{R}^2;\ 0<x<\alpha(t),\ t\in(0,T)
	\}.
	$$

	We investigate a multi-objective controllability problem concerning the linear one-dimensional wave equation
	%
	\begin{equation}\label{WE1}
		\left\{
		\begin{array}{lll}
			u_{tt} - u_{xx} = 0 & \mbox{in} & Q_{T,\alpha},\\
			u(\cdot,0) = u_0,\quad u_t(\cdot,0)= u_1 &\mbox{in}& (0,1),
		\end{array}
		\right.
	\end{equation}
	where $u_0, u_1$ are initial data and $\alpha$ is a $C^2(\mathbb{R}_+)$ function with the following properties:
	\begin{enumerate}[label=(H\arabic*)]
		\item $\alpha(0)=1$;
		\item $\alpha'(t)\in (m,M)$, $\forall t\in \mathbb{R}_+$, with $(m, M) \subset (0, 1)$;
		\item $\alpha'$ is monotone.
	\end{enumerate}
	
	The lateral boundary $\Sigma_{T,\alpha}$ of $Q_{T,\alpha}$ is composed by two disjoint sets $\Sigma_{T}^0$ and $\Sigma_{T}^\alpha$, given by
	$$
	\Sigma_{T}^0 = \{(0,t);\ t\in (0,T)\}
	\mbox{\quad and\quad }
	\Sigma_{T}^\alpha = \{(\alpha(t),t);\ t\in (0,T)\},
	$$
	that is,  $\Sigma_{T,\alpha}=\Sigma_{T}^0\cup\Sigma_{T}^\alpha$. We are interested in studying the scenario where the solution $u$ of the wave equation \eqref{WE1} is subjected to either one of the following boundary conditions:
	%
		%
		%
	%
	\begin{equation}\label{BC0}
		u = \left\{\begin{array}{lll}
			\displaystyle f+v & \mbox{on} & \Sigma_{T}^0,
			\\\noalign{\smallskip}
			\displaystyle 0 & \mbox{on} & \Sigma_{T}^\alpha,
		\end{array}
		\right.
	\end{equation}
	or
	%
		%
		%
		%
		%
		\begin{equation}\label{BCalpha}
			u = \left\{\begin{array}{lll}
				\displaystyle 0 & \mbox{on} & \Sigma_{T}^0,
				\\\noalign{\smallskip}
				\displaystyle g+w & \mbox{on} & \Sigma_{T}^\alpha,
			\end{array}
			\right.
		\end{equation}
		where $f,v \in L^2(\Sigma_{T}^{0})$, $g, w \in L^2(\Sigma_{T}^{\alpha})$ are control functions.
		
		In the aforementioned context, we then outline the problems to be addressed in this study.
		
		For each control, $f$ and $v$, acting within the wave system \eqref{WE1}, \eqref{BC0}, and for $g$ and $w$ in \eqref{WE1}, \eqref{BCalpha}, we define specific objectives:
		
		\begin{enumerate}[label=]
			\item {\bf \underline{Task for $f$ and $g$}:} Approximate controllability
			
			The objective is to determine optimal controls $f$ and $g$ that effectively approximate the solution $u$ of equations \eqref{WE1}, \eqref{BC0}, and \eqref{WE1}, \eqref{BCalpha}, respectively, to a predefined target $u^T$ at the final time $T$. In other words, the aim is to ensure that the system satisfied by $u$ achieves approximate controllability.
			
			\item {\bf \underline{Task for $v$ and $w$}:} Minimization
			
			The objective is to achieve the optimal approximation (as measured by distance calculated using norms) of the solution $u$ of equations \eqref{WE1}, \eqref{BC0} to another predefined datum $u_2$, utilizing the smallest possible selection of $v$. This scenario can be reformulated as the minimization of an appropriate functional. Similarly, for a solution $u$ of equations \eqref{WE1}, \eqref{BCalpha}, the objective remains unchanged.
		\end{enumerate}
		
		Our problem consist in to analyze the existence of pairs $f,v$ that accomplish the tasks above simultaneously. Since each task has a different objective to be achieved by $u$, then finding the pair $f,v$ is called a multi-objective problem.
		The same approach will be considered for the system \eqref{WE1}, \eqref{BCalpha}, where the tasks of $f$ and $v$ are given to $g$ and $w$, respectively.
		
		

		Let us consider the functionals $\mathcal{J}: L^2(\Sigma_T^0)\rightarrow \mathbb{R}$ and $\mathcal{K}: L^2(\Sigma_{T}^{\alpha})\rightarrow \mathbb{R}$ given by
		\begin{equation}
			\mathcal{J}(f) = \dfrac{1}{2}\int_{\Sigma_T^0} |f|^2\,d\Sigma,
			\quad\quad
			\mathcal{K}(g) = \dfrac{1}{2}\int_{\Sigma_{T}^{\alpha}} |g|^2\,d\Sigma,
		\end{equation}
		and the functionals $\mathcal{J}_f: L^2(\Sigma_T^0)\rightarrow \mathbb{R}$ and $\mathcal{K}_g: L^2(\Sigma_{T}^{\alpha})\rightarrow \mathbb{R}$ given by
		\begin{equation}
			\mathcal{J}_f (v) = \dfrac{1}{2}\displaystyle\int_{Q_{T,\alpha}} |u(x,t; f,v, u_0,u_1) - u_2|^2\,dx\,dt + \dfrac{\tilde\sigma}{2}\int_{\Sigma_T^0} |v|^2\,d\Sigma,
		\end{equation}
		\begin{equation}
			\mathcal{K}_g (w) = \dfrac{1}{2}\displaystyle\int_{Q_{T,\alpha}} |u(x,t; g,w, u_0,u_1) - u_4|^2\,dx\,dt + \dfrac{\tilde\mu}{2}\int_{\Sigma_{T}^{\alpha}} |w|^2\,d\Sigma.
		\end{equation}
		In the definition of the four functionals above, $\tilde\sigma,\tilde\mu$ are positive constants to be specified later; $u(x,t; f,v,u_0,u_1)$ is the solution of system \eqref{WE1}, \eqref{BC0} associated to the data $f,v,u_0,u_1$; and $u(x,t; g,w,u_0,u_1)$ stands for the solution of  system \eqref{WE1}, \eqref{BCalpha} associated to the data  $g,w,u_0,u_1$.

		Our main theoretical results are the following:
		
		\begin{theorem}\label{M1}
			Let $u_0\in L^2(0,1)$, $u_1\in H^{-1}(0,1)$,  $u^T\in L^2(0,\alpha(T))$, and $u_2\in L^2(Q_{T,\alpha})$. Let us assume that conditions $(H1)-(H3)$ hold true and
			\begin{equation}\label{T0}
				T > T_1 = \dfrac{1}{M}\left(e^{\tfrac{2M^2(1-m)}{m(1-M)^3}} - 1\right).
			\end{equation}
			Then, if $\tilde\sigma$ is large enough, for each $\varepsilon>0$ there exist controls $f_\varepsilon\in L^2(\Sigma_T^0)$ and $v_\varepsilon\in L^2(\Sigma_T^0)$ such that
			\begin{enumerate}
				\item $\|u(\cdot,T;f_\varepsilon,v_\varepsilon,u_0,u_1) - u^T\|_{L^2(0,\alpha(T))}<\varepsilon$;
				\item $\displaystyle\mathcal{J}_{f_\varepsilon} (v_\varepsilon)=\min_{L^2(\Sigma_T^0)} \mathcal{J}_{f_\varepsilon} (v)$;
				\item $\displaystyle\mathcal{J}(f_\varepsilon)=\min_{L^2(\Sigma_T^0)} \mathcal{J}(f)$.
			\end{enumerate}
		\end{theorem}
		
		\begin{theorem}\label{M2}
			Let $u_0\in L^2(0,1)$, $u_1\in H^{-1}(0,1)$, $u^T\in L^2(0,\alpha(T))$, and $u_4\in L^2(Q_{T,\alpha})$. Let us assume that conditions $(H1)-(H3)$ hold true and
			$$
			T > T_2 = \dfrac{1}{M}\left(e^{\tfrac{2M^2(1-m)(1+M)}{m(1-M)^2}} - 1\right).
			$$
			Then, if $\tilde\mu$ is large enough, for each $\varepsilon>0$ there exist controls $g_\varepsilon\in L^2(\Sigma_T^{\alpha})$ and $w_\varepsilon\in L^2(\Sigma_T^{\alpha})$ such that
			\begin{enumerate}
				\item $\|u(\cdot,T;g_\varepsilon,w_\varepsilon,u_0,u_1) - u^T\|_{L^2(0,\alpha(T))}<\varepsilon$;
				\item $\displaystyle\mathcal{K}_{g_\varepsilon} (w_\varepsilon)=\min_{\Sigma_T^{\alpha}} \mathcal{K}_{g_\varepsilon}(w)$;
				\item $\displaystyle\mathcal{K}(g_\varepsilon)=\min_{\Sigma_T^{\alpha}} \mathcal{K}(g)$.
			\end{enumerate}
		\end{theorem}

		In Theorems \ref{M1} and \ref{M2}, conditions $1-3$ delineate the objectives that must be fulfilled by the pairs of controls $(f,v)$ and $(g,w)$. Condition $1$ signifies the well-established property of approximate controllability.
		
		In order to solve our multi-objective problem concerning the system \eqref{WE1}, \eqref{BC0} we will assume a cooperative scenario between the controls $f$ and $v$ where $v$ depends on $f$. This choice of scenario configures the called Stackelberg hierarchical strategy (see \cite{stalck}) and, due to the framework, we usually name $f$ as the leader control and $v$ as the follower control. This hierarchical assignment of the controls allows us to proceed as follows. Firstly, for any fixed $f$ we prove the existence of a control $v$ uniquely defined by $f$ and that minimize $\mathcal{J}_f$. Because of such unique writing of $v=v(f)$, we streamline the search of $(f,v)$ into finding only a control $f$. This is made by using an auxiliary system that provides a characterization of $v$. The mentioned auxiliary system, combined with \eqref{WE1}, \eqref{BC0}, generates the so-called optimality system, that has only the control $f$ acting. Therefore, proving the approximate controllability for the optimality system we find controls satisfying the conditions $1-2$. To finalize, we minimize the functional $\mathcal{J}$ in the set composed by all controls of the optimality system. This final step furnish the pair of controls that solve Theorem \ref{M1}. The same approach will be used to prove Theorem \ref{M2}.
		
		In addition to the Stackelberg hierarchical cooperative strategy, there are several other approaches available for solving multi-objective problems involving partial differential equations. For instance, we can consider the non-cooperative optimization strategy proposed by Nash \cite{nash} and the Pareto cooperative strategy \cite{pareto}, both of which stem from game theory and are primarily motivated by economics. In the realm of partial differential equations, the concept of hierarchical control was introduced in \cite{lions1}, where only one leader and one follower were considered.

		
		In the context of approximate controllability, several multi-objective problems were studied. We can cite the extension of  \cite{DL} in D\'iaz \cite{D1}, that uses Fenchel-Rockafellar duality theory to give
 a characterization of the solution; the papers \cite{RA1, RA2} of Ramos et al. that studied Nash equilibrium for linear parabolic PDEs and for the Burger's equation; and Guill\'en-Gonz\'alez et al. in \cite{GO} that developed the Stackelberg-Nash approximate controllability for the Stokes system. In the context of hyperbolic equations, Lions \cite{lions2} applied Stackelberg techniques from \cite{stalck} and solved an approximate controllability problem for a linear wave equation. In this approach, were considered one primary controller (the leader) and additional secondary control (the follower). Conversely, the hierarchical exact controllability properties for hyperbolic PDEs have been extensively explored in recent years; for example, see \cite{Luc} and \cite{Luc1}.

		
		
		A noteworthy contribution presented here is to provide an improvement to the results presented in \cite{A,B}. There, the author dealt with the problem of approximate controllability for
 a hyperbolic equation in the one-dimensional case using a very particular boundary function, namely, $\alpha(t) = 1+kt.$ Here we intend to extend his results to the case of a much more general boundary
 function with some hypotheses. More specifically, this work covers functions much more comprehensive than straight lines, for example, $\alpha(t) = 1 + (t + \arctan t)/c$, where $c$ is any constant
  greater than 2. In the case of $k = 1,$ some results have been obtained in \cite {Cui1}. On the other hand in the case  $k > 1,$ the moving boundary is a spacelike surface, on which an initial condition rather than a boundary condition needs to be imposed.  
		
		This article is structured as follows: In Section \ref{WP}, we provide a concise overview of the well-posedness of the problems. Section \ref{Char} is dedicated to reorganizing the problems with the goal of deriving the optimality systems necessary for proving Theorems \ref{M1} and \ref{M2}. Sections \ref{Th1} and \ref{Th2} are dedicated to proving Theorem \ref{M1} and Theorem \ref{M2}, respectively.
		The final section is reserved for comments an conclusions.

		\section{Well-posedness of the problem}\label{WP}
		
		The objective of this section is to establish the well-posedness of the problems represented by equations (\ref{WE1}), (\ref{BC0}) and (\ref{WE1}), (\ref{BCalpha}). To begin with, we perform a change of variables to convert the problems represented by equations \eqref{WE1}, \eqref{BC0} and \eqref{WE1}, \eqref{BCalpha} from a non-cylindrical domain into a cylindrical domain. Upon isolating one of the lateral boundaries of the cylindrical domain in each problem, we encounter four controls, namely $f$, $v$, $g$, and $w$, acting on the separation. Here, we consider $f$ and $g$ as leader controls, while $v$ and $w$ serve as followers, following the terminology of Stackelberg's strategy.
		
		Let us observe that when $\displaystyle (x,t)$ varies in $Q_{T,\alpha}$, the point $\displaystyle (y,t) = (x/\alpha(t), t)$, varies in $\displaystyle Q=\Omega \times (0,T),$ where
		$\Omega = (0,1)$. Then, considering the diffeomorphism $\displaystyle \zeta:Q_{T, \alpha} \to Q$ defined by $\zeta(x,t)=(y,t)$, 
		we set $\Gamma = \zeta(\Sigma), \Gamma_0 = \zeta(\Sigma_{T}^0), \Gamma_{\alpha} = \zeta(\Sigma_{T}^{\alpha}).$
		%
		%
		Therefore, the change of variables $\displaystyle u(x,t)=z(y,t)$, transforms the boundary value problems \eqref{WE1}, \eqref{BC0} and \eqref{WE1}, \eqref{BCalpha} into the equivalent systems
		\begin{equation}\label{eq1.14}
			\begin{array}{lll}
				\displaystyle z_{tt} +  L z =  0 & \mbox{in} & Q,
				\\\noalign{\smallskip}
				\displaystyle z(y,t) = (f+v)\mathbf{1}_{\Gamma_0}
				&\mbox{on}& \Gamma,
				\\\noalign{\smallskip}
				\displaystyle z(y, 0) = z_0,\quad  z_{t}(y,0) = z_1 & \mbox{in} & \Omega,
			\end{array}
		\end{equation}
		and
		%
		\begin{equation}\label{eq1.15}
			\begin{array}{lll}
				\displaystyle z_{tt} +  L z =  0 & \mbox{in} & Q,
				\\\noalign{\smallskip}
				\displaystyle z(y,t) = (g+w)\mathbf{1}_{\Gamma_{\alpha}}
				&\mbox{on}& \Gamma,
				\\\noalign{\smallskip}
				\displaystyle z(y, 0) = z_0,\quad  z_{t}(y,0) = z_1 & \mbox{in} & \Omega,
			\end{array}
		\end{equation}
		respectively, where
		\begin{equation}\label{coe}
			\begin{array}{lll}
				\displaystyle Lz = - \Big[\frac{\beta(y,t)}{\alpha(t)}z_{y} \Big]_{y} + \frac{\gamma(y,t)}{\alpha(t)}z_{ty} + \frac{\tau(y,t)}{\alpha(t)}z_y, \\[10pt]
				\displaystyle \beta(y,t) = \frac{1  - \alpha'^2(t)y^2}{\alpha(t)}, \quad \gamma(y,t) = -2\alpha'(t)y, \quad \tau(y,t) = -\alpha''(t)y ,\\[10pt]
				\displaystyle z_{0}(y) = u_{0}(x), \quad z_{1}(y) = u_{1}(x) + \alpha'(0)yu_{x}(0) , \\[10pt]
				\displaystyle \Gamma = \Gamma_0 \cup \Gamma_{\alpha}, \; {\Gamma}_0 = \{(0,t) : t \in (0, T) \}, \; {\Gamma}_{\alpha} = \{(1,t) : t \in (0, T) \},
			\end{array}
		\end{equation}
		with regularity
		\begin{equation}\label{reg}
			\beta \in C^1(\overline{Q}), \quad \gamma, \tau \in W^{1, \infty}(\Omega).
		\end{equation}
		
		This change of variables also defines the functionals.
		
		\begin{equation}
			{J}(f) = \dfrac{1}{2}\int_{\Gamma_0} |f|^2\,d\Gamma,
			\quad\quad
			{K}(g) = \dfrac{1}{2}\int_{\Gamma_{\alpha}} |g|^2\,d\Gamma,
		\end{equation}
		and
		\begin{equation}\label{K2}
			{J}_{f}(v) = \dfrac{1}{2}\displaystyle\int_{Q} \alpha(t) |z(y, t; f, v, z_0, z_1) - z_2|^2\,dy\,dt + \dfrac{\sigma}{2}\int_{\Gamma_0} |v|^2\,d\Gamma,
		\end{equation}
		
		\begin{equation}\label{J2}
			{K}_{g}(w) = \dfrac{1}{2}\displaystyle\int_{Q} \alpha(t) |z(y, t; g, w, z_0, z_1) - z_4|^2\,dy\,dt + \dfrac{\mu}{2}\int_{\Gamma_{\alpha}} |w|^2\,d\Gamma,
		\end{equation}
		where $\sigma, \mu > 0$ are constants and $z_2(y,t), z_4(y, t)$ are given functions in $L^2(Q).$
		
		
		\begin{remark}
			For each set of data $z_0 \in L^2(\Omega)$, $z_1 \in H^{-1}(\Omega)$, and controls $f$, $v \in L^2(\Gamma_0)$, $g$, $w \in L^2(\Gamma_{\alpha})$, there exists a unique solution bu transposition $z$ of the problems $(\ref{eq1.14})$ and $(\ref{eq1.15})$, as discussed in \cite{Mii}.
			This solution has the regularity
			\begin{equation*}
				\displaystyle z \in C\big( [0,T]
				;L^2(\Omega)) \cap C^{1}\big( [0,T] ; H^{-1}(\Omega)).
			\end{equation*}
			
		\end{remark}

		Hence, it follows that the problems represented by equations (\ref{WE1}), (\ref{BC0}) and (\ref{WE1}), (\ref{BCalpha}) are also well-posed. Consequently, the functionals $\mathcal{J}_f$, $\mathcal{K}_g$, $J_f$, and $K_g$ are also well-posed.
		

		\section{Caracterization of the problems}\label{Char}
		In this section, our objective is to derive the optimality system, which characterizes the follower controls and will serve as the foundation for proving the main theorems. Given any fixed leader controls $f \in L^2(\Gamma_0)$ and $g \in L^2(\Gamma_{\alpha})$, we aim to establish the existence and uniqueness of solutions to the problems.
		\begin{equation} \label{eq3.10}
			\begin{array}{l}
				\displaystyle\inf_{\widehat{v} \in L^2(\Gamma_0)}{J}_f(\widehat{v})
			\end{array}
		\end{equation}
		and
		\begin{equation} \label{eq3.11}
			\begin{array}{l}
				\displaystyle\inf_{\widehat{w} \in L^2(\Gamma_{\alpha})}{K}_g(\widehat{w}),
			\end{array}
		\end{equation}
		and also a characterization of such solutions using an adjoint system.
		
		
		
		Since, for each fixed $f \in L^2(\Gamma_0)$ and $g\in L^2(\Gamma_{\alpha})$, the functionals ${J}_f$, ${K}_g$ are $C^1$, coercive, weakly lower semicontinuous and strictly convex, there exists a unique solution $v \in L^2({\Gamma}_0)$ to \eqref{eq3.10} and $w \in L^2({\Gamma}_{\alpha})$ to \eqref{eq3.11}, i.e.,
		\begin{equation*}
			\displaystyle {J}_f(v)= \inf_{\widehat{v} \in L^2({\Gamma}_0)} {J}_f(\widehat{v}), \quad \displaystyle {K}_g(w)= \inf_{\widehat{w} \in L^2({\Gamma}_{\alpha})} {K}_g(\widehat{w}).
		\end{equation*}
		In particular, ${J}_{f}'(v) = {K}_{g}'(w) = 0$, i.e,
		\begin{equation} \label{eq3.20}
			\int_{0}^{T}\int_{\Omega}{\alpha(t)}(z-z_2)\widehat{z}dy\,dt + \sigma\int_{\Gamma_0}v\widehat{v}d\Gamma = 0, \;\;\forall\, \widehat{v} \in L^2(\Gamma_0),
		\end{equation}
		where $\widehat{z}$ is solution of
		%
		\begin{equation}\label{eq3.21}
			\begin{array}{lll}
				\displaystyle \widehat z_{tt} +L\widehat{z} = 0 & \mbox{in} & Q,
				\\\noalign{\smallskip}
				\displaystyle \widehat z(y,t) = \widehat{v}\mathbf{1}_{\Gamma_{0}},
				&\mbox{on}& \Gamma,
				\\\noalign{\smallskip}
				\displaystyle \widehat z(y, 0) = 0,\quad  \widehat z_{t}(y,0) = 0 & \mbox{in} & \Omega,
			\end{array}
		\end{equation}
		and
		\begin{equation*} \label{eq3.22}
			\int_{0}^{T}\int_{\Omega}{\alpha(t)}(z-z_4)\tilde{z}dy\,dt + \mu\int_{\Gamma_{\alpha}}w\widehat{w}d\Gamma = 0, \;\;\forall\, \widehat{w} \in L^2(\Gamma_{\alpha}),
		\end{equation*}
		where $\tilde{z}$ is solution of
		%
		\begin{equation}\label{eq3.23}
			\begin{array}{lll}
				\displaystyle \tilde z_{tt} +L\tilde{z} = 0 & \mbox{in} & Q,
				\\\noalign{\smallskip}
				\displaystyle \tilde z(y,t) = \widehat{w}\mathbf{1}_{\Gamma_{\alpha}}
				&\mbox{on}& \Gamma,
				\\\noalign{\smallskip}
				\displaystyle \tilde z(y, 0) = 0,\quad  \tilde z_{t}(y,0) = 0 & \mbox{in} & \Omega.
			\end{array}
		\end{equation}
		
		Let us consider the adjoint system to \eqref{eq3.21} and \eqref{eq3.23} defined by
		%
		\begin{equation}\label{sac}
			\begin{array}{lll}
				\displaystyle p_{tt} + L^{\ast}\,p = \alpha(t)\left(z - z_2\right) & \mbox{in} & Q,
				\\\noalign{\smallskip}
				\displaystyle p(y, t) = 0
				&\mbox{on}& \Gamma
				\\\noalign{\smallskip}
				\displaystyle p(y, T) = p_{t}(y, T) = 0 & \mbox{in} & \Omega,
			\end{array}
		\end{equation}
		where
\begin{equation*}
	\begin{array}{l}
		\displaystyle L^{\ast}\,p =  -\left[\dfrac{\beta(y, t)}{\alpha(t)} p_y \right]_y + \left[\dfrac{\gamma(y, t)}{\alpha(t)} p\right]_{yt} - \left[\dfrac{\tau(y, t)}{\alpha(t)} p\right]_y
	\end{array}
\end{equation*}
is the formal adjoint of the operator $\displaystyle L$, and $\beta(y,t),$ $\gamma(y,t),$ and $\tau(y,t)$ are given in (\ref{coe}). \\ Multiplying $\eqref{sac}_1$ by $\widehat{z}$ and integrating by parts, we get that
\begin{equation*}
	\int_{0}^{T}\int_\Omega \alpha(t)(z - z_2)\widehat{z}\,dy\,dt + \int_{\Gamma} \frac{\beta(y, t)}{\alpha(t)}\,p_y\,\widehat{z}\,d\Gamma = 0.
\end{equation*}
Separating the boundary $\Gamma = \Gamma_0 \cup \Gamma_{\alpha}$ in the integral and using \eqref{eq3.21}, we have
\begin{equation*} \label{eq3.33}
	\int_{0}^{T}\int_\Omega \alpha(t)(z - z_2)\widehat{z}\,dy\,dt + \int_{\Gamma_0} \frac{1}{\alpha^2(t)}\,p_y\,\widehat{v}\,d\Gamma = 0,
\end{equation*}
so that by \eqref{eq3.20} we get
\begin{equation*}\label{ci}
	v = \dfrac{1}{\sigma \alpha^2(t)} p_y \ \ \mbox{ on } \ \ \Gamma_0.
\end{equation*}
Similarly, we can find an explicit expression for $w$, i.e,
\begin{equation*}\label{cii}
	w = \dfrac{1}{\mu \alpha^2(t)} p_y \ \ \mbox{ on } \ \ \Gamma_{\alpha}.
\end{equation*}
Therefore, we have proven the following theorem.
\begin{theorem}[Optimality system for the follower control] \label{teN} For each $\displaystyle f \in L^2(\Gamma_0)$ and $\displaystyle g \in L^2(\Gamma_{\alpha})$ there exists a unique Nash equilibrium $\displaystyle v$ and $w$ in the sense of \eqref{eq3.10} and \eqref{eq3.11}, respectively. Moreover, the followers $\displaystyle v, w$ are given by
	\begin{equation}\label{cseg}
		\displaystyle \displaystyle v = v(f)=\frac{1}{\sigma \alpha^2(t)}\,\;p_y\;\;\mbox{ on }\;\;\Gamma_0,
	\end{equation}
	
	\begin{equation}\label{csegg}
		\displaystyle \displaystyle w = w(g)=\frac{1}{\mu \alpha^2(t)}\,\;p_y\;\;\mbox{ on }\;\;\Gamma_{\alpha},
	\end{equation}
	where in the first case $\displaystyle \{ z,p \}$ is the unique solution of
	%
	\begin{equation}\label{eq3.40}
		\begin{array}{lll}
			\displaystyle z_{tt} + Lz = 0, \quad p_{tt} + L^{\ast}\,p = \alpha(t)\left(z - z_2\right) & \mbox{in} & Q,
			\\\noalign{\smallskip}
			\displaystyle z = \left(f + \frac{1}{\sigma \alpha^2(t)}\;\,p_y\right)\mathbf{1}_{\Gamma_{0}}, \quad p = 0
			&\mbox{on}& \Gamma,
			\\\noalign{\smallskip}
			\displaystyle z(0) = z_0,\quad z_{t}(0) = z_1, \quad  p(T) = p_{t}(T) = 0 & \mbox{in} & \Omega,
			\\\noalign{\smallskip}
		\end{array}
	\end{equation}
	and in the second case $\displaystyle \{ z,p \}$ is the unique solution of
	%
	\begin{equation}\label{eq3.38}
		\begin{array}{lll}
			\displaystyle z_{tt} + Lz = 0, \quad p'' + L^{\ast}\,p = \alpha(t)\left(z - z_4\right) & \mbox{in} & Q,
			\\\noalign{\smallskip}
			\displaystyle z = \left(g + \frac{1}{\tilde \sigma \alpha^2(t)}\;\,p_y\right)\mathbf{1}_{\Gamma_{\alpha}}, \quad p = 0
			&\mbox{on}& \Gamma,
			\\\noalign{\smallskip}
			\displaystyle z(0) = z_0,\quad z_{t}(0) = z_1, \quad p(T) = p'(T) = 0 & \mbox{in} & \Omega.
			\\\noalign{\smallskip}
		\end{array}
	\end{equation}
	
\end{theorem}

\section{Proof of the Theorem \eqref{M1}}\label{Th1}
After characterizing the follower control, our focus shifts to proving the existence of a leader control $f$ as a solution to the problem
\begin{equation} \label{inf11}
	\displaystyle \inf_{f \in \mathcal{U}(v^0, v^1, \varepsilon)} {J}(f),
\end{equation}
where $\displaystyle \mathcal{U}(v^0, v^1, \varepsilon)$, is the set of admissible controls
\begin{equation*}\label{admcon}
	\displaystyle \mathcal{U}(v^0, v^1, \varepsilon)=\{f \in L^2({\Gamma_0}); \; \|\{z(\cdot, T; f, v), z_t(\cdot, T; f, v)\}-\{v^0, v^1\}\|_{L^2(\Omega)\times H^{-1}(\Omega)} < \varepsilon\},
\end{equation*}
with $v$ as in \eqref{cseg}, $v^0$ and $v^1$ are arbitrary in $L^2(\Omega) \times H^{-1}(\Omega)$.

Since \eqref{eq1.14} is linear, it suffices to prove this for null initial data. We begin by decomposing the solution $\displaystyle (z,p)$ of \eqref{eq3.40} as
\begin{equation} \label{eq3.39}
	(z, p) = (\nu_0, p_0) + (h, q)
\end{equation}
where ($\nu_0$, $p_0$) is the solution of
\begin{equation}\label{eq3.46}
	\begin{array}{lll}
		\displaystyle (\nu_0)_{tt} + L\,\nu_0 = 0, \quad (p_0)_{tt} + L^{\ast}p_0 = \alpha(t) \left(\nu_0 - z_2\right) & \mbox{in} & Q,
		\\\noalign{\smallskip}
		\displaystyle \nu_0 = \frac{1}{\sigma \alpha^2(t)}\,({p_0})_{y}\mathbf{1}_{\Gamma_{0}}, \quad p_0 = 0
		&\mbox{on}& \Gamma,
		\\\noalign{\smallskip}
		\displaystyle \nu_0(0) = (\nu_0)_t(0) = 0, \quad p_0(T) = (p_0)_{t}(T) = 0 & \mbox{in} & \Omega,
	\end{array}
\end{equation}
and $\displaystyle (h, q)$ is the solution of
\begin{equation}\label{eq3.48}
	\begin{array}{lll}
		\displaystyle h_{tt} + L\,h = 0, \quad q_{tt} + L^{\ast}q = \alpha(t) h & \mbox{in} & Q,
		\\\noalign{\smallskip}
		\displaystyle h = \left(f + \frac{1}{\sigma \alpha^2(t)}\,q_{y}\right)\mathbf{1}_{\Gamma_{0}}, \quad q = 0
		&\mbox{on}& \Gamma,
		\\\noalign{\smallskip}
		\displaystyle h(0) = h_t(0) = 0, \quad q(T) = q_t(T) = 0 & \mbox{in} & \Omega.
	\end{array}
\end{equation}
Let $A$ be the continuous linear operator defined by
\begin{equation} \label{eq3.44}
	\begin{array}{ccll}
		A \ : & \! L^2(\Gamma_0) & \! \longrightarrow & \! H^{-1}(\Omega) \times L^2(\Omega) \\
		& \! f & \! \longmapsto & \! Af = \big\{ h_t(T; f), -h(T; f) \big\},
	\end{array}
\end{equation}
Consider the convex proper functions $F_1:L^2(\Gamma_0) \longrightarrow \mathbb{R} \cup \{\infty\}$ and $F_2: H^{-1}(\Omega) \times L^2(\Omega) \longrightarrow \mathbb{R} \cup \{\infty\}$ defined by
\begin{equation*} \label{eq3.119}
	\displaystyle F_1(f) = \frac{1}{2} \int_{\Gamma_0}|f|^{2}\,d\Gamma,
\end{equation*}
and
\begin{equation*} \label{eq3.120}
	\begin{array}{ccl}
		F_2(Af) = F_2\big(\{h_t(T, f), -h(T, f)\}\big) \\
		= \left\{
		\begin{array}{l}
			0, \text{ if }
			\left\{
			\begin{array}{l}
				h_t(T) \in - (\nu_0)_{t}(T) + B_{H^{-1}(\Omega)}(v^1, \varepsilon),\\
				-h(T) \in \nu_0(T) - B_{L^2(\Omega)}(v^0, \varepsilon),
			\end{array}
			\right.\\
			+ \infty, \text{ otherwise}.
		\end{array}
		\right.
	\end{array}
\end{equation*}
Then, problem \eqref{inf11} becomes equivalent to
\begin{equation} \label{eq3.122}
	\begin{array}{l}
		\displaystyle \inf_{f \in L^2(\Gamma_0)}\big[F_1(f) + F_2(Af)\big],
	\end{array}
\end{equation}
once we establish that the range of $A$ is dense in $H^{-1}(\Omega) \times L^2(\Omega)$.

To establish this, we employ a density criterion from \cite{aubin}, that is, if the following holds:
\begin{equation*}
	\displaystyle \Big\langle \big\langle A f , \{f^0, f^1\} \big\rangle \Big\rangle = 0 \quad \forall f \in L^2(\Gamma_0) \Rightarrow \{f^0, f^1\} = \{0, 0\},
\end{equation*}
where $\Big\langle \big\langle . , . \big\rangle \Big\rangle$ represent the duality
pairing between $ H^{-1}(\Omega) \times L^2(\Omega) $ and $ H_{0}^{1}(\Omega) \times
L^2(\Omega)$.

Given $\{f^0, f^1\} \in H_0^1(\Omega) \times L^2(\Omega)$, let us consider the adjoint system of \eqref{eq3.48}:
\begin{equation} \label{eq3.52}
	\begin{array}{lll}
		\displaystyle \varphi_{tt} +  L^\ast \varphi = \alpha(t) \psi, \quad \psi_{tt} +  L \, \psi = 0 & \mbox{in} & Q,
		\\\noalign{\smallskip}
		\displaystyle \varphi = 0, \quad \psi = \frac{1}{\sigma \alpha^2(t)}\,\varphi_{y} \mathbf{1}_{\Gamma_0}
		\quad
		&\mbox{on}& \Gamma,
		\\\noalign{\smallskip}
		\displaystyle \varphi(T) = f^0,\quad  \varphi_{t}(T) = f^1, \quad \psi(0) = \psi_{t}(0) = 0 & \mbox{in} & \Omega.
	\end{array}
\end{equation}
Multiplying by $q$ the equation  satisfied by $\psi$ in \eqref{eq3.52}$_1$, where $q$
solves (\ref{eq3.48}), and integrating over $Q$, we obtain
\begin{equation}\label{eq3.99}
	\int_{0}^{T}\int_{\Omega} \psi q_{tt}\,dy\,dt + \int_{0}^{T}\int_{\Omega} \psi L^{*}q\,dy\,dt + \int_{\Gamma} \dfrac{\beta(y, t)}{\alpha(t)} q_y \psi \, d\Gamma = 0.
\end{equation}
Using the equation satisfied by $q$ in \eqref{eq3.48}$_1$ and the boundary condition of $\psi$ on $\Gamma$ in \eqref{eq3.52}$_2$, then \eqref{eq3.99} becomes:
\begin{equation} \label{eq3.47}
	\int_{0}^{T}\int_{\Omega} \alpha(t) h\,\psi\,dy\,dt =- \frac{1}{\sigma}\int_{\Gamma_0} \frac{1}{\alpha^4(t)}\,q_{y}\, \varphi_{y} d\Gamma.
\end{equation}
Multiplying by $h$ the equation satisfied by $\varphi$ in \eqref{eq3.52}, where $h$
solves \eqref{eq3.48}, and integrating over $Q$, we
obtain:
\begin{equation*}
	\int_0^T\int_{\Omega} \alpha(t) \psi h \, dy\,dt = \big( h(T),f^1 \big) - \langle h_t(T),f^0 \rangle_{H^{-1}(\Omega) \times  H_{0}^{1}(\Omega)} - \int_{\Gamma} \dfrac{\beta(y, t)}{\alpha(t)} \varphi_y h \, d\Gamma = 0.
\end{equation*}
Using the condition of $h$ on $\Gamma$ in \eqref{eq3.48}$_2$, we get
\begin{equation}\label{eq28}
	\begin{array}{ll}
		\displaystyle \int_0^T\int_{\Omega} \alpha(t) \psi h \, dy\,dt = \big( h(T),f^1 \big) - \langle h_t(T),f^0 \rangle_{H^{-1}(\Omega) \times  H_{0}^{1}(\Omega)}\\ \displaystyle - \int_{\Gamma_0}\dfrac{1}{\alpha(t)^2} \, \varphi_y \, f\,d\Gamma - \displaystyle \dfrac{1}{\sigma} \int_{\Gamma_0}\dfrac{1}{\alpha(t)^4} \, \varphi_y \, q_y\, d\Gamma.
	\end{array}
\end{equation}
Combining the equations \eqref{eq3.47} and \eqref{eq28}, we can conclude:
\begin{equation}\label{eq3.49}
	-\int_{\Gamma_0} \dfrac{1}{\alpha^2(t)} \varphi_y f d\,\Gamma = \langle h_t(T), f^0 \rangle - (h(T), f^1).
\end{equation}

Considering the right-hand side of equation \eqref{eq3.49} as the inner product between $ H^{-1}(\Omega) \times L^2(\Omega) $ and $
H_{0}^{1}(\Omega) \times L^2(\Omega)$ of $\displaystyle \{h_t(T), -h(T)\}$ with $\{f^0, f^1\}$, we obtain:
\begin{equation*}
	\Big\langle \big\langle A f , \{f^0, f^1\} \big\rangle \Big\rangle = - \int_{\Gamma_0}\frac{1}{\alpha^2(t)}\, \varphi_{y}\,f\,d\Gamma.
\end{equation*}
Thus, if
\begin{equation*}
	\Big\langle \big\langle A f , \{f^0, f^1\} \big\rangle \Big\rangle = 0 \quad \forall f \in L^2(\Gamma_0),
\end{equation*}
then $\varphi_{y}= 0 \mbox{ on } \Gamma_0,$ which implies $\psi \equiv 0$.

Therefore, $\varphi$ satisfies
\begin{equation*}\label{eq3.54}
	\begin{array}{lll}
		\displaystyle \varphi_{tt} +  L^\ast \varphi =  0 & \mbox{in} & Q,
		\\\noalign{\smallskip}
		\displaystyle \varphi(y,t) = 0,
		\quad
		\displaystyle \varphi_{y}\mathbf{1}_{\Gamma_0}(y, t) = 0
		&\mbox{on}& \Gamma,
		\\\noalign{\smallskip}
		\displaystyle \varphi(y,T) = f^0,\quad \varphi_t(y,T) = f^1 & \mbox{in} & \Omega.
	\end{array}
\end{equation*}
According to Theorem 3.6 of \cite {Cui12},
we conclude that $\{\widehat{f}^0, \widehat{f}^1\} = \{0, 0\}$ if $T > T_1$. This concludes the proof of density.

We now resume the proof of the problem
\eqref{eq3.122}. Using the Fenchel-Rockafellar Theorem \cite{R}(see also \cite{Bre}, \cite{EK}), we obtain
\begin{equation} \label{eq3.124}
	\begin{array}{l}
		\displaystyle \inf_{f \in L^2(\Gamma_0)}[F_1(f) + F_2(Af)]\\[5pt]\displaystyle  = -\inf_{\{f^0, f^1\} \in H_{0}^{1}(\Omega) \times L^2(\Omega)} [F_{1}^{*}\big(A^*\{f^0, f^1\}\big) + F_{2}^{*}\{-f^0, -f^1\}],
	\end{array}
\end{equation}
where $\displaystyle A^*: H_{0}^{1}(\Omega) \times L^2(\Omega) \longrightarrow L^2(\Gamma_0)$ is the adjoint of $\displaystyle A$, and $\displaystyle F_i^*$ is the convex conjugate function of $\displaystyle F_i$ for $i = 1, 2$. These are given by
\begin{equation*} \label{eq3.121}
	A^*\{f^0, f^1\} = -\dfrac{1}{\alpha^2(t)}\,\varphi_{y}, \quad F_{1}^{*}(f) = F_1(f),	
\end{equation*}
where $\varphi$ is the solution of (\ref{eq3.52}), and
\begin{equation*}\label{eq3.125.2}
	\begin{array}{ccl}
		\displaystyle  F_{2}^{*}(\{f^0, f^1 \})  &=& \displaystyle \langle v^1 - (\nu_0)_{t}(T), f^0\rangle_{H^{-1}(\Omega) \times H_{0}^{1}(\Omega)} \\&&+\big( \nu_0(T) - v^0 , f^1\big)
		\displaystyle + \varepsilon\|f^0\| + \varepsilon|f^1|.
	\end{array}
\end{equation*}
Then, we can rewrite \eqref {eq3.124} as
\begin{equation} \label{eq3.1240}
	\begin{array}{l}
		\displaystyle \inf_{f \in L^2(\Gamma_0)}[F_1(f) + F_2(Af)] = -\inf_{\{f^0, f^1\} \in  H_{0}^{1}(\Omega) \times L^2(\Omega)}  \Theta \big(\{f^0, f^1\}\big),
	\end{array}
\end{equation}
where the functional $\Theta : H_{0}^{1}(\Omega) \times L^2(\Omega) \longrightarrow \mathbb{R}$ is given by
\begin{equation*}
	\begin{array}{l}
		\displaystyle \Theta \big(\{f^0, f^1\}\big) = \frac{1}{2}\int_{\Gamma_0}\left(\frac{1}{\alpha^2(t)}\right)^2 \varphi_{y}^2d\;\Gamma  + \big( v^0 - \nu_0(T), f^1\big) \\
		- \langle v^1 - (\nu_0)_{t}(T), f^0\rangle_{H^{-1}(\Omega) \times H_{0}^{1}(\Omega)} + \varepsilon\|f^0\| + \varepsilon|f^1|.
	\end{array}
\end{equation*}
Since $\Theta$ is a continuous, coercive, and strictly convex functional, there exists a unique solution $\{f^0, f^1\}$ of
\begin{equation*}
	\inf_{\{f^0,f^1\} \in  H_{0}^{1}(\Omega) \times L^2(\Omega)}  \Theta \big(\{f^0, f^1\}\big).
\end{equation*}
Then $\{f^0, f^1\}$ satisfies the variational inequality
\begin{equation}\label{eq801}
	\begin{array}{l}
		\displaystyle \int_{\Gamma_0} \left(\frac{1}{\alpha^2(t)}\right)^2 \varphi_y (\widehat{\varphi}_y - \varphi_y) d\,\Gamma  +  \left\langle\left\langle \{-\eta^1, \eta^0\}, \{\widehat{f}^0, \widehat{f}^1\} \right\rangle\right\rangle + \varepsilon\|\widehat{f}^0\|\\
		+ \varepsilon|\widehat{f}^1| - \left\langle\left\langle \{-\eta^1, \eta^0\}, \{f^0, f^1\} \right\rangle\right\rangle - \varepsilon\|f^0\| - \varepsilon|f^1| \geq 0,
	\end{array}
\end{equation}
for all $\{\widehat{f}^0, \widehat{f}^1\} \in H_0^1(\Omega) \times L^2(\Omega)$.
From \eqref{eq3.49} with $f = - \frac{1}{\alpha^2(t)} \varphi_y$, we have
\begin{equation}\label{eq199}
	\begin{array}{l}
		\displaystyle \langle h_t(T),f^0 \rangle_{H^{-1}(\Omega) \times  H_{0}^{1}(\Omega)} - \big( h(T),f^1 \big) = \displaystyle \int_{\Gamma_0}\left(\frac{1}{\alpha^2(t)}\right)^2\,\varphi_{y}^2\,d\Gamma,
	\end{array}
\end{equation}
futhermore,
\begin{equation}\label{eq200}
	\begin{array}{l}
		\displaystyle \langle h_t(T),\widehat{f}^0 \rangle_{H^{-1}(\Omega) \times  H_{0}^{1}(\Omega)} - \big( h(T),\widehat{f}^1 \big) = \displaystyle \int_{\Gamma_0}\left(\frac{1}{\alpha^2(t)}\right)^2\,\varphi_{y} \widehat{\varphi}_y\,d\Gamma.
	\end{array}
\end{equation}
Combining \eqref{eq199} and \eqref{eq200} in \eqref{eq801}, we get
\begin{equation*}
	\begin{array}{l}
		\displaystyle \langle z_t(T) - v^1, \widehat{f}^0 - f^0 \rangle - (z(T) - v^0, \widehat{f}^1 - f^1))\\
		+ \varepsilon(\|\widehat{f}^0\| - \|f^0\|) + \varepsilon(|\widehat{f}^1| - |f^1|) \geq 0, \quad \forall \{\widehat{f}^0, \widehat{f}^1\} \in H_0^1(\Omega) \times L^2(\Omega).
	\end{array}
\end{equation*}
The variational inequality above implies that $f = - \frac{1}{\alpha^2(t)} \varphi_y$ is the desired function in \eqref{inf11}. Thus, the proof is concluded.

\section{Proof of the Theorem \eqref{M2}}\label{Th2}
The proof of Theorem \eqref{M2} follows a similar structure to that of Theorem \eqref{M1}. Our objective is to determine a leader control $g$ as the solution of
$g$ as the solution of
\begin{equation} \label{eq3.7cil.1}
	\displaystyle \inf_{g \in \mathcal{W}(v^0, v^1, \varepsilon)} {K}(g),
\end{equation}
where $\displaystyle \mathcal{W}(v^0, v^1, \varepsilon)$ is the set of admissible controls
\begin{equation*}\label{admconn}
	\displaystyle \mathcal{W}(v^0, v^1, \varepsilon)=\{g \in L^2({\Gamma_{\alpha}}); \; \|\{z(\cdot, T; g, w), z_t(\cdot, T; g, w)\}-\{v^0, v^1\}\|_{L^2(\Omega)\times H^{-1}(\Omega)} < \varepsilon\},
\end{equation*}
with $w$ as in \eqref{csegg}. Since \eqref{eq1.15} is linear, it suffices to prove this for null initial data.
	We start decomposing the solution $\displaystyle (z,p)$ of \eqref{eq3.38} as
	\begin{equation*} \label{eq5.39}
		(z, p) = (\tilde \nu_0, \tilde p_0) + (\tilde h, \tilde q)
	\end{equation*}
	where $\tilde \nu_0$, $\tilde p_0$ are given as solutions of
	\begin{equation*}
		\begin{array}{lll}
			\displaystyle \tilde \nu_0'' + L\,\tilde \nu_0 = 0, \quad \tilde p_0'' + L^{\ast}\tilde p_0 = \alpha(t) \left(\tilde \nu_0 - z_4\right) & \mbox{in} & Q,
			\\\noalign{\smallskip}
			\displaystyle \tilde\nu_0 = \frac{1}{\tilde\sigma \alpha^2(t)}\,({\tilde p_0})_{y}\mathbf{1}_{\Gamma_{\alpha}}, \quad \tilde p_0 = 0
			&\mbox{on}& \Gamma,
			\\\noalign{\smallskip}
			\displaystyle \tilde \nu_0(0) = \tilde \nu_0'(0) = 0, \quad \tilde p_0(T) = \tilde p_0'(T) = 0 & \mbox{in} & \Omega,
		\end{array}
	\end{equation*}
	and $\displaystyle \tilde h, \tilde q$ are given as solutions of
	\begin{equation*}\label{eq3.488}
		\begin{array}{lll}
			\displaystyle \tilde h'' + L\,\tilde h = 0, \quad \tilde q'' + L^{\ast}\tilde q = \alpha(t) \tilde h & \mbox{in} & Q,
			\\\noalign{\smallskip}
			\displaystyle \tilde h = \left(g + \frac{1}{\tilde\sigma \alpha^2(t)}\,\tilde q_{y}\right)\mathbf{1}_{\Gamma_{\alpha}}, \quad \tilde q = 0
			&\mbox{on}& \Gamma,
			\\\noalign{\smallskip}
			\displaystyle \tilde h(0) = \tilde h'(0) = 0, \quad \tilde q(T) = \tilde q'(T) = 0 & \mbox{in} & \Omega.
		\end{array}
	\end{equation*}
	We set
	\begin{equation*} \label{eq5.44}
		\begin{array}{ccll}
			A \ : & \! L^2(\Gamma_{\alpha}) & \! \longrightarrow & \! H^{-1}(\Omega) \times L^2(\Omega) \\
			& \! g & \! \longmapsto & \! Ag = \big\{ \tilde h'(T; g), -\tilde h(T; g) \big\},
		\end{array}
	\end{equation*}
	which defines
	\begin{equation*}
		A \in \mathcal{L}\left( L^2(\Gamma_{\alpha}); \;H^{-1}(\Omega) \times L^2(\Omega)\right).
	\end{equation*}
	
	Consider the convex proper functions $G_1:L^2(\Gamma_{\alpha}) \longrightarrow \mathbb{R} \cup \{\infty\}$ and $G_2: H^{-1}(\Omega) \times L^2(\Omega) \longrightarrow \mathbb{R} \cup \{\infty\}$ defined by
	\begin{equation*}
		\displaystyle G_1(g) = \frac{1}{2} \int_{\Gamma_{\alpha}}|f|^{2}\,d\Gamma,
	\end{equation*}
	and
	\begin{equation*}
		\begin{array}{ccl}
			G_2(Ag) = F_2\big(\{\tilde h'(T, g), -\tilde h(T, g)\}\big) \\
			= \left\{
			\begin{array}{l}
				0, \text{ if }
				\left\{
				\begin{array}{l}
					\tilde h'(T) \in - \tilde \nu_0'(T) + B_{H^{-1}(\Omega)}(v^1, \varepsilon),\\
					-\tilde h(T) \in \tilde \nu_0(T) - B_{L^2(\Omega)}(v^0, \varepsilon),
				\end{array}
				\right.\\
				+ \infty, \text{ otherwise}.
			\end{array}
			\right.
		\end{array}
	\end{equation*}
	Then, problem \eqref{eq3.7cil.1} becomes equivalent to
	\begin{equation} \label{eq3.1222}
		\begin{array}{l}
			\displaystyle \inf_{g \in L^2(\Gamma_{\alpha})}\big[G_1(g) + G_2(Ag)\big],
		\end{array}
	\end{equation}
	once we establish that the range of $A$ is dense in $H^{-1}(\Omega) \times L^2(\Omega)$. The proof of this density follows the same approach as in the previous proof, and in the conclusion, we utilize Theorem 3.5 from \cite{Cui12}.

	Returning to the proof of \eqref{eq3.1222}, we use the Fenchel-Rockafellar Theorem to obtain
	\begin{equation*}
		\begin{array}{l}
			\displaystyle \inf_{g \in L^2(\Gamma_{\alpha})}[G_1(g) + G_2(Ag)] = -\inf_{\{g^0, g^1\} \in  H_{0}^{1}(\Omega) \times L^2(\Omega)}  \Theta \big(\{g^0, g^1\}\big),
		\end{array}
	\end{equation*}
	where the functional $\Theta : H_{0}^{1}(\Omega) \times L^2(\Omega) \longrightarrow \mathbb{R}$ is given by
	\begin{equation*}
		\begin{array}{l}
			\displaystyle \Theta \big(\{g^0, g^1\}\big) = \frac{1}{2}\int_{\Gamma_{\alpha}}\left(\frac{1}{\alpha^2(t)}\right)^2 \varphi_{y}^2d\;\Gamma  + \big( v^0 - \tilde \nu_0(T), g^1\big) \\
			- \langle v^1 - \tilde \nu_0'(T), g^0\rangle_{H^{-1}(\Omega) \times H_{0}^{1}(\Omega)} + \varepsilon\|g^0\| + \varepsilon|g^1|.
		\end{array}
	\end{equation*}
	The solution $\{g^0, g^1\}$ of
	\begin{equation*}
		\inf_{\{g^0,g^1\} \in  H_{0}^{1}(\Omega) \times L^2(\Omega)}  \Theta \big(\{g^0, g^1\}\big)
	\end{equation*}
	satisfies the variational inequality
	\begin{equation*}
		\begin{array}{l}
			\displaystyle \langle z'(T) - v^1, \widehat{g}^0 - g^0 \rangle - (z(T) - v^0, \widehat{g}^1 - g^1))\\
			+ \varepsilon(\|\widehat{g}^0\| - \|g^0\|) + \varepsilon(|\widehat{g}^1| - |g^1|) \geq 0, \quad \forall \{\widehat{g}^0, \widehat{g}^1\} \in H_0^1(\Omega) \times L^2(\Omega),
		\end{array}
	\end{equation*}
	with $g = - \frac{1}{\alpha^2(t)} \varphi_y$.
	
	The variational inequality above implies that $g = - \frac{1}{\alpha^2(t)} \varphi_y$ is the desired function in \eqref{eq3.7cil.1}. Thus, the proof is concluded.

\section{Additional comments and conclusions}\label{Conclusion}
In this paper, we have demonstrated the feasibility of achieving approximate controllability results for the wave equation through the implementation of a Stackelberg strategy. The theoretical motivation behind our study is to enhance the findings presented in \cite{A,B}. Specifically, this work encompasses a broader range of functions compared to simple straight lines, exemplified by functions like $\alpha(t) = 1 + (t + \arctan t)/c$, where $c$ is any constant greater than 2. Importantly, the concepts and techniques elucidated in this paper can be extrapolated to various analogous scenarios, such as hierarchical approximate controllability for linear and semilinear heat and plate equations, diverse forms of non-cylindrical control domains, and analogous boundary control problems.

Concurrently, recent theoretical advancements have been made in controllability utilizing Stackelberg-Nash strategies (refer to \cite{Luc2},\cite{Luc}). Exploring the extension of these concepts to other multi-objective control problems, including the computation of Stackelberg-Pareto equilibria, presents an intriguing avenue for future research.

\paragraph{\bf Acknowledgement}

 Pedro Paulo A.  Oliveira is grateful for the support received from CAPES-Brazil. This research constitutes a component of the Ph.D. thesis of Pedro Paulo A. Oliveira at the Department of Mathematics of
  the Federal University of Piau\'{i}. The author extends appreciation to the host institution for its warm hospitality. Isa\'{i}as P. de Jesus acknowledges partial support from grants provided by
  CNPq/Brazil [Grants: 307488/2019-5, 305394/2022-3] and PRPG (UFPI).





%
\paragraph{\bf Declarations}
\paragraph {Conflicts of interest}  On behalf of all authors the corresponding author states that there is no conflict of interest.

\paragraph

\end{document}